\documentclass[12pt]{amsart}
\usepackage{amscd,amssymb,graphicx}
\author{Alice Fialowski}
\address{E\"otv\"os Lor\'and University\\
Budapest, Hungary} \email{fialowsk@cs.elte.hu}
\author{Michael Penkava}
\address{University of Wisconsin\\
Eau Claire, WI 54702-4004} \email{penkavmr@uwec.edu}
\thanks{The research of the authors was partially supported by
 OTKA grants T043641 and T043034 and by grants from the University of Wisconsin-Eau Claire.}%

\newtheorem{thm}{Theorem}[section]

\newtheorem{lma}[thm]{Lemma}
\theoremstyle{definition}

\def \ph{\varphi}

\def \refeq#1{equation (\ref{#1})}
\def \ra{\rightarrow}

\def \hom{\mbox{\rm Hom}}

\def \tns{\otimes}

\def \mcom{,\cdots,}
\def \k{\mbox{$\mathbb K$}}

\def \C{\mbox{$\mathbb C$}}
\def \Z{\mbox{$\mathbb Z$}}

\def\zt{\mbox{$\Z_2$}}

\def\sh{\operatorname{Sh}}

\def\ad{\operatorname{ad}}

\def\inv{^{-1}}
\def\d{d}

\def\im{\operatorname{Im}}

\def\dl{\delta}

\def\linf{\mbox{$L_\infty$}}
\def\and{\mbox{ \rm and }}

\def\s#1{(-1)^{#1}}

\def\pha#1#2{\ph^{#1}_{#2}}

\def\dl#1#2{\delta^{#1}_{#2}}
\def\psa#1#2{\psi^{#1}_{#2}}

\def\thmref#1{Theorem (\ref{#1})}
\def\inv{^{-1}}

\def\dbar#1{\overline{\overline{#1}}}
\def\dl{\delta+\lambda}
\def\bdl{\bar\delta+\bar\lambda}
\def\bbdl{\overline{\overline{\delta+\lambda}}}
\def\hdlm{H_{\mu,\delta+\lambda}}
\def\hdlmp{H_{\mu,\delta+\lambda,\psi}}
\def\hdlpm{H_{\psi,\delta+\lambda,\mu}}
\def\hrest{H_\mu(\ker(D_{\dl}))}
\def\hmu{H_\mu}
\def\hpsi{H_\psi}
\def\hdlp{H_{\psi,\delta+\lambda}}
\def\cbdl#1{\left[\overline{#1}\right]}
\def\cbdlp#1{\left[\overline{\overline{#1}}\right]}
\def\ctpsi#1{\left\{\cbdl{#1}\right\}}
\def\ctmu#1{\left\{\cbdlp{#1}\right\}}
\def\grest{\mbox{$G_{\text{rest}}$}}
\def\ggen{\mbox{$G_{\text{gen}}$}}
\def\gdiag{\mbox{$G_\Delta$}}
\def\gdiagmd{\mbox{$\gdiag(\mu,\delta)$}}
\def\ggenmd{\mbox{$\ggen(\mu,\delta)$}}
\def\gdiagmdl{\mbox{$\gdiag(\mu,\delta,\lambda)$}}
\def\ggenmdl{\mbox{$\ggen(\mu,\delta,\lambda)$}}
\def\GL{{\mathbf{GL}}}

\begin{document}
\setlength{\multlinegap}{0pt}
\title{Extensions  of (super) Lie algebras}%

\address{}%
\email{}%

\thanks{}%

\date{\today}
\begin{abstract}
In this paper, we give a purely cohomological interpretation of the
extension problem for (super) Lie algebras; that is the problem of
extending a Lie algebra by another Lie algebra. We then give a
similar interpretation of infinitesimal deformations of extensions.
In particular, we consider infinitesimal deformations of representations
of a Lie algebra.
\end{abstract}
\maketitle

\section{Introduction}
The notion of an extension of a (super) Lie algebra is an old idea
in Lie theory (see \cite{ce,Lev-Sen,mac}), even though not much
explicit research has been done recently (see
\cite{amr,neeb1,neeb2}). One of the difficulties may have been the
lack of tools which help to identify and classify extensions of a
given algebra. As we will point out in this paper, to give a
satisfactory definition of the moduli space of all nonequivalent
extensions of a Lie algebra, it is necessary to consider a more
general notion of equivalence of extensions than usually appears in
the literature.

The structure of this paper is as follows. After some preliminary
definitions and explanation of notation in section \ref{prelim}, we
will recall the classical definition of a super Lie algebra
extension  and discuss the notion of equivalent extensions in
section \ref{extensions}.  In section \ref{genext}, we introduce a
more general concept of equivalence of extensions, and give the
classification of extensions of a Lie algebra in purely
cohomological terms. To indicate the advantage of the general
definition of equivalence, we give a few examples which could not be
obtained with the restricted definition. In section \ref{defext}, we
define infinitesimal deformations of extensions of super Lie
algebras, and in section \ref{infdef} we introduce infinitesimal
deformations of representations.

Our main guide in studying the extension problem comes from the
recent studies of the moduli spaces of equivalence classes of Lie
algebra structures on low dimensional spaces in
\cite{fp3,fp8,Ott-Pen}. The set of equivalence classes of extensions
of a Lie algebra by another Lie algebra determines a moduli space as
well, and the notion of equivalence we use is based on the
philosophy that the moduli space of extensions should fit nicely
into the moduli space of Lie algebra structures on the underlying
vector space given by the sum of the two Lie algebras. A too
restricted notion of equivalence leads to a moduli space which is
too large.

\section{Preliminaries}\label{prelim}
In classical Lie theory, the cohomology of a Lie algebra is studied by
considering a differential on the space of maps from the \emph{exterior algebra}
of the vector space to the vector space itself; that is, the \emph{cochains with
values in the adjoint representation}. The exterior algebra $\bigwedge V$ of a non-graded
vector space $V$ has a natural $\Z$-grading, which induces a \zt-grading on the exterior algebra.
There is also a superbracket, the Richardson-Nijenhuis bracket \cite{nr},
on this space of cochains, equipping it with the
structure of a super Lie algebra. Thus,
even for ordinary Lie algebras, one needs to introduce the notion of a \zt-graded
space in order to understand deformation theory.

In this paper, we will find it is more convenient to
consider Lie algebra structures as codifferentials in
the space of coderivations of a \emph{symmetric coalgebra},
rather than the exterior algebra, a language which is particularly
well suited for generalization to other contexts, but which may not be familiar
to Lie algebra theorists.  Accordingly, we will provide here an explanation for
our terminology. In fact, in the case of ordinary (not super) Lie algebras, only
a dictionary need be provided, because the codifferentials we consider are simply
Lie algebra structures, the coderivations we consider are simply elements of the cochain
complex of the Lie algebra, and the \zt-grading of the coderivations is given only in
terms of the degree of the cochain.  The bracket of coderivations coincides with the
usual Richardson-Nijenhuis bracket of cochains.

For super Lie algebras, it is also possible to study the cohomology without introducing
the coderivation point of view, but the bracket of cochains is much simpler to express
in terms of the bracket of coderivations; the signs which appear are simpler to express and easier
to understand. The fact that the bracket structure of the cochains has a natural interpretation
as a bracket of coderivations of a symmetric coalgebra was first noticed by J. Stasheff
in \cite{sta4}.

Let $V$ be a \zt-graded vector space over a field $\k$, which, in our examples,
we will assume is $\C$, although, for the most part, $\k$ can be any field
whose characteristic is not 2 or 3. Recall that a super Lie algebra structure on
$V$ is given by a graded antisymmetric map $V\tns V\ra V$, denoted by
$(a,b)\mapsto [a,b]$, which satisfies the graded Jacobi identity
\begin{equation*}
[a,[b,c]]=[[a,b],c]+\s{ab}[b,[a,c]],
\end{equation*}
for homogeneous elements. If $V=V_e$, then the structure is a Lie algebra, and
the sign $\s{ab}$ does not appear.

Recall that for a Lie algebra, the cochain complex with values in the adjoint
representation, which is important for deformation theory, is given by
$C(V)=\prod_{n=0}^\infty C^n(V)$, where
$C^n(V)=\hom(\bigwedge^n V,V)$.
The space $C(V)$ has a
\zt-grading, where the parity of an element in $C^n$ is even if $n$ is odd,
and odd if $n$ is even. There is a bracket operation on $C(V)$, which satisfies
$[C^k(V),C^l(V)]\subseteq C^{k+l-1}(V)$, and this bracket equips $C(V)$ with the structure
of a \zt-graded Lie algebra.
In fact, if $\ph\in C^k(V)$, $\psi\in C^l(V)$, and $n=k+l-1$
 then we can define $\phi\circ\psi\in
C^{n}(V)$ by
\begin{equation*}
\phi\circ\psi(v_1\mcom v_{n})=\!\!
\!\!\!\sum_{\sigma\in\sh(l,k-1)}\s{\sigma}
\phi(\psi(v_{\sigma(1)}\mcom v_{\sigma(l)}),v_{\sigma(l+1)}\mcom v_{\sigma(n)}),
\end{equation*}
where $\s{\sigma}$ is the sign of the permutation $\sigma$ and $\sh(l,k-1)$ are the
permutations in $\Sigma_n$ which are increasing on the first $l$ and the last $k-1$
elements.  The bracket of cochains is given by
\begin{equation*}
[\phi,\psi]=\phi\circ\psi -\s{\phi\psi}\psi\circ\phi,
\end{equation*}
where $\s{\phi\psi}=\s{(k-1)(l-1)}$.

The cohomology of the Lie algebra is determined by
the differential $D:C(V)\ra C(V)$, given by $D(\phi)=[l,\phi]$, where $l$ is the Lie
algebra structure, interpreted as an element of $C^2(V)$. The fact that $D^2=0$ follows
from the Jacobi identity, which in terms of cochains is the codifferential property
$[l,l]=0$. Note that $l$ is an odd cochain, because it lies in $C^2(V)$.  Therefore,
by the graded Jacobi identity on $C(V)$, we have
\begin{equation*}
D^2(\phi)=[l,[l,\phi]]=[[l,l],\phi]-[l,[l,\phi]=-[l,[l,\phi]].
\end{equation*}
The formula above explains why we do not allow a field of characteristic 2 in our theory.
The reason fields of characteristic 3 are a problem is more subtle, but has to do with
the fact that we want the triple bracket $[\alpha,[\alpha,\alpha]]$ to vanish for any
$\alpha$, and this will be true as long as the characteristic of the field is not 3.
Even when the Lie algebra is not graded, the cochain complex $C(V)$ is a \zt-graded
Lie algebra, so fields of characteristic 2 or 3 are still problematic.

When working with \zt-graded Lie algebras, it is convenient to consider the
\emph{parity reversion} $W=\Pi V$, where $W_o=V_e$ and $W_e=V_o$. We can
identify $C(V)=\hom(\bigwedge V,V)$ with $C(W)=\hom(S(W),W)$, where $S(W)$ is
the \zt-graded symmetric coalgebra of $W$. Note that the symmetric coalgebra coincides
with the symmetric algebra as a vector space, equipped with the
cocommutative, coassociative coproduct $\Delta:S(W)\ra S(W)\tns S(W)$, given by
\begin{equation*}
\Delta(w_1\cdots w_n)=\sum_{\substack{k=0\dots n\\\sigma\in\sh(k,n-k)}}\epsilon(\sigma)
w_1\cdots w_k\tns w_{k+1}\cdots w_n,
\end{equation*}
where $\epsilon(\sigma)$ is a sign determined by the rule
\begin{equation*}
w_{\sigma(1)}\cdots w_{\sigma(n)}=\epsilon(\sigma)w_1\cdots w_n.
\end{equation*}

There is a natural identification of $C(W)$ with the coderivations of $S(W)$, and the
bracket in $C(W)$ corresponds to the bracket of coderivations. In fact, if
$\phi\in C^k(W)$, $\psi\in C^l(W)$, and $n=k+l-1$, then we define $\phi\circ\psi\in C^n(W)$
by
\begin{equation*}
\phi\circ\psi(w1\cdots w_n)=
\sum_{\sigma\in\sh(l,k-1)}
\epsilon(\sigma)\phi(\psi(w_{\sigma(1)}\cdots w_{\sigma(l)})w_{\sigma(l+1)}\cdots w_n),
\end{equation*}

Then
\begin{equation*}
[\phi,\psi]=\phi\circ\psi-\s{\phi\psi}\psi\circ\phi,
\end{equation*}
where $\s{\phi\psi}$ is a sign determined by the product of the parities of $\phi$ and $\psi$, rather than their
degrees as cochains.

Note that when $V$ is an ordinary Lie algebra, $\epsilon(\sigma)=\s{\sigma}$ is just the
sign of the permutation, and the parity of the element in $C(W)$ associated to an element
in $C^k(V)$ is just $k-1$, so the bracket in $C(W)$ is just the bracket in $C(V)$. Therefore,
for ordinary Lie algebras, the construction of $C(W)$ is just a change of language.

If $d$ is the element in $C^2(W)$ associated
to a Lie algebra structure $l$ in $C^2(V)$, then $d$ satisfies the following
version of the Jacobi identity.
\begin{equation*}
d(d(a,b),c)+\s{bc}d(d(a,c),b)+\s{a(b+c)}d(d(b,c),a)=0.
\end{equation*}
In fact, $d$ is an odd coderivation of the symmetric coalgebra $S(W)$, and moreover the
Jacobi identity above is simply the codifferential property
$[d,d]=0$. We say that $d$ is an odd codifferential on $S(W)$.

One defines the cohomology of the codifferential $d$ on $W$ in terms of the
differential $D:C(W)\ra C(W)$ given by $D(\ph)=[d,\ph]$. We call the resulting homology
the cohomology $H(d)$ of the super Lie algebra structure on $V$ given by the codifferential
$d$ on $W$.

The main advantage in working with $C(W)$ instead of $C(V)$ is that in $C(W)$, one
only has to consider the parity of an element, instead of its parity and its degree
as a cochain. Thus we only have to consider the structure of $C(W)$ as a graded Lie
algebra, which is a real advantage in computations.

Let $\{w_1\mcom w_n\}$ be a basis of the space $W$, and let $w_{i,j}=w_iw_j$ denote
the product of two basis elements, which lies in $S^2(W)$. Consider the map
$$S(i):\{1\mcom n\}\ra\{1\mcom n\}\times\{1\mcom n\}$$ given by the sequence of ordered pairs
$$(1,2),(1,3),(2,3),(1,4),(2,4)\mcom (n-1,n).$$
Then there is a natural basis
$\{w_{S(1)}\mcom w_{S{\binom n2}}\}$ for $S^2(W)$, which we will
use in order to represent codifferentials $d:S^2(W)\ra W$ as matrices.
Moreover, any element of $C^2(W)$ can be also represented as a sum
$d=a_i\psa{S(i)}j$, where
$\psa{i,j}k:S^2(W)\ra W$ is defined by $\psa{i,j}k(w_{l,m})=\delta^{i,j}_{l,m}w_k$,
and we shall also express our codifferentials using this second notation
in our examples.

\section{Extensions of Lie Algebras}\label{extensions}
We consider the extension picture of a (super) Lie algebra $W$
by a Lie algebra $M$. In other words we consider the short exact sequence
of Lie algebras
\begin{equation}\label{exactseq}
0\ra M\ra V\ra W\ra 0,
\end{equation}
For convenience, we introduce the following notation for certain spaces of cochains on $V=M\oplus W$.
\begin{align*}
C^n=&\hom(W^n,W)\\
C^{k,l}=&\hom(M^kW^l,M),
\end{align*}
where $M^kW^l$ is the subspace of $S^{k+l}(V)$ determined by
products of $k$ elements from $M$ and $l$ elements from $W$.
The Lie
algebra structure on $V$ is determined by the following maps:
\begin{align*}
&\delta\in C^2=\hom(W^2,W):& \quad \text{the Lie algebra structure on $W$}\\
&\psi\in C^{0,2}=\hom(W^2,M):& \quad \text{the 2-cocycle on $M$}\\
&\lambda\in C^{1,1}=\hom(MW,M):& \quad \text{the $W$-``module
structure'' on $M$}\\
&\mu\in C^{2,0}=\hom(M^2,M):& \quad \text{the Lie algebra structure
on $M$}
\end{align*}
The fact that $d$ has no terms from $\hom(MW\oplus M^2,W)$
reflects the fact that $M$ is an ideal in $V$.
The Jacobi identity is the condition that the
coderivation $d=\delta+\lambda+\rho+\mu+\psi$ is a codifferential on
$S(V)$, which means that $[d,d]=0$. Now, in general, we see that
$[d,d]\in\hom(S^3(V),V)$. By decomposing this space and
considering which parts the brackets of the terms $\delta$, $\lambda$
$\mu$ and $\psi$ are defined on, we obtain
\begin{align}
\label{lrel1} &[\delta,\delta]=0:
\text{ The Jacobi identity for $\delta$.}\\
\label{lrel2} &[\mu,\mu]=0:
\text{ The Jacobi identity for $\mu$.}\\
\label{lrel3} &[\delta,\lambda]+1/2[\lambda,\lambda]+[\mu,\psi]=0:
\text{ The ``module'' relation.}\\
\label{lrel4} &[\mu,\lambda]=0:
\text{ The module-algebra structures are compatible.}\\
\label{lrel5} &[\delta+\lambda,\psi]=0:
\text{ $\psi$ is a 2-cocycle with values in $M$.}
\end{align}
We also have the automatic relations
\begin{math}
[\mu,\delta]=[\psi,\psi]=0.
\end{math}
When $[\mu,\psi]=0$, \refeq{lrel3} says that $\lambda$ is a module
structure in the usual sense, so this relation can be considered as
saying that $\lambda$ determines a generalized module structure.
 Note that $\delta+\mu$ is a codifferential, the direct sum of the
codifferentials on $M$ and $W$, and $\lambda+\psi$ satisfies the
Maurer-Cartan formula (MC-formula)
\begin{equation}
[\delta+\mu,\lambda+\psi]+\tfrac12\left[\lambda+\psi,\lambda+\psi\right]=0,
\end{equation} with respect to
the codifferential $\delta+\mu$.  Furthermore, $\mu+\lambda+\psi$ satisfies the MC formula

\begin{equation}
[\delta,\mu+\lambda+\psi]+\tfrac12[\mu+\lambda+\psi,\mu+\lambda+\psi]=0,
\end{equation}
with respect to the codifferential $\delta$.
All of these basic facts about Lie algebra
extensions are well known.
(See, for example \cite{ce,mac}.)
We summarize the main facts in the theorem below.
\begin{thm}\label{th1}
If $\delta$ is a Lie algebra structure on $W$ and $\mu$ is a Lie algebra structure on $M$,
then $d=\delta+\mu+\lambda+\psi$, where $\lambda\in\hom(MW,M)$ and $\psi\in\hom(W^2,M)$ determines
a Lie algebra structure on $M\oplus W$ precisely when the three conditions below hold:
\begin{align}
&[\delta,\lambda]+\frac12[\lambda,\lambda]+[\mu,\psi]=0\label{mc-twist}\\
&[\mu,\lambda]=0\label{mu-cocycle}\\
&[\delta+\lambda,\psi]=0\label{ker-cond}
\end{align}
\end{thm}
\section{Equivalence of extensions of Lie algebras}\label{genext}
A (restricted) equivalence of extensions of  Lie algebras is given by a commutative diagram
of the form
$$
\begin{CD}
0@>>>M@>>>V@>>>W@>>>0\\
@.@|@VVfV@|\\
0@>>>M@>>>V@>>>W@>>>0\\
\end{CD}
$$
where we assume that in  the top row, $V$ is equipped with the codifferential
$d=\delta+\mu+\lambda+\psi$, and in the bottom row, it is equipped with the codifferential
$d'=\delta+\mu+\lambda'+\psi'$, and $f$ is a morphism of Lie algebras (which is necessarily
an isomorphism). In order for the diagram to commute, we must have $f(m,w)=(m+\beta(w),w)$,
where $\beta:W\ra M$ is linear. Note  that if $g(m,w)=(m-\beta(w),w)$, then $g=f\inv$.

The requirement that $f$ be a morphism of Lie algebras is simply
that $d'=g^*d=g\inv\circ d\circ g$, where $g$  is extended to  an
automorphism of $S(V)$. It follows that
\begin{align*}
\delta'=&\delta\\
\mu'=&\mu\\
\lambda'=&\lambda-[\mu,\beta]\\
\psi'=&\psi-[\delta+\lambda-\tfrac12[\mu,\beta],\beta].
\end{align*}
Since $d'=g^*d$, it  is
automatically true that $d'$ is a codifferential.

If we consider the map $\beta$ as an endomorphism of $V$, then $\beta^2=0$, and therefore
$\exp(\beta)=1+\beta=f$. But then $f^*$ is given by the formula
$$
f^*=\exp(-\ad_\beta),
$$
which means that the action of $f^*$ on $\hom(S(V),V)$ is expressible in terms of
brackets.  Note that
\begin{align*}
\ad_\beta:C^k\ra C^{0,k}\\
\ad_\beta:C^{k,l}\ra C^{k-1,l+1}.
\end{align*}
The term $\tfrac12[[\mu,\beta],\beta]$ is simply
$\tfrac12\ad_\beta^2(\mu)$, arising from the exponential, and this
is the only nonvanishing  term in $f^*(d)$ involving a double
bracket with $\beta$, and all higher order brackets vanish.

An extension $V$ is said to be a semidirect product provided that $W$ is a subalgebra of
$V$. This just means that $\psi=0$. Clearly, an extension is equivalent to a semidirect
product precisely when there is some linear  $\beta:W\ra M$ such that
$
\psi=[\delta+\lambda-\tfrac12[\mu,\beta],\beta]
$.
Let us summarize the results above in the following theorem, which is given in a slightly different form in
\cite{amr}.
\begin{thm}\label{th2}

If $d=\delta+\mu+\lambda+\psi$ and $d'=\delta+\mu+\lambda'+\psi'$
are two extensions of the Lie algebra structure $\delta$ on $W$ by a Lie algebra
structure $\mu$ on $M$, then they are equivalent extensions precisely when there
is a $\beta\in\hom(W,M)$ such that
\begin{align}\label{equivext}
\lambda'=&\lambda-[\mu,\beta]\\
\psi'=&\psi-[\delta+\lambda-\tfrac12[\mu,\beta],\beta].
\end{align}
Moreover the extension is equivalent to a semidirect product precisely when
there is some $\beta\in\hom(W,M)$ such that
\begin{equation}
\psi=[\delta+\lambda-\tfrac12[\mu,\beta],\beta].
\end{equation}
\end{thm}

The group $\grest$ of automorphisms $f$ of $V$ of the form
$$f(m,w)=(m+\beta(w),w)$$ acts on the set of extensions. The set of
equivalence classes of such extensions under this group action is
the \emph{restricted moduli space of extensions of $\delta$ by
$\mu$}.
\subsection{Infinitesimal extensions and infinitesimal equivalence}
The notion of an infinitesimal extension is expressed in the form
$$d=\delta+\mu+t(\lambda+\psi),$$ where $t^2=0$ is an infinitesimal parameter. The conditions for
$d$ to be an extension reduce to
\begin{align*}
&[\delta,\lambda]+[\mu,\psi]=0\\
&[\mu,\lambda]=0\\
&[\delta,\psi]=0
\end{align*}
If $\alpha\in C(V)$, then denote $D_\alpha=\ad_\alpha$. When
$\alpha$ is odd and $[\alpha,\alpha]=0$, then $D^2_\alpha=0$, so
$D_\alpha$ is called a \emph{coboundary operator} on $C(V)$, and
$H_\alpha=\ker(D_\alpha)/\im(D_\alpha)$ is the \emph{cohomology}
induced by $\alpha$. An element $\phi$ such that $D_\alpha(\phi)=0$
is called a \emph{$D_\alpha$-cocycle}, and its image $[\phi]_\alpha$
in $H_\alpha$ is a \emph{$D_\alpha$-cohomology class}. An element of
the form $D(\phi)$ is called a $D_\mu$-coboundary. The bracket on
$C(V)$ descends to a bracket on $H_\alpha$, so $H_\alpha$ inherits
the structure of a Lie superalgebra.  Since $\delta$, $\mu$ and
$\psi$ are all codifferentials, they determine coboundary operators.

Note that there is a certain symmetry in the roles of $\delta$ and
$\mu$ in the above conditions for an infinitesimal extension, in the
sense that if you interchange the codifferentials $\delta$ and
$\mu$, and the cochains $\lambda$ and $\psi$, then the conditions
are unchanged. We have
\begin{align*}
D_\delta:C^k\ra C^{k+1}&\quad&D_\delta:C^{k,l}\ra C^{k,l+1}\\
D_\mu:C^k\ra 0&\quad& D_\mu:C^{k,l}\ra C^{k+1,l}.
\end{align*}
Since $[\delta,\mu]=0$, it follows that $D_\delta$ and $D_\mu$
anticommute. As a consequence,
\begin{equation*}
D_\mu:\ker(D_\delta)\ra \ker(D_\delta),
\end{equation*}
so we can define the cohomology $H_\mu(\ker\delta)$ determined by
the restriction of $D_\mu$ to $\ker(D_\delta)$.

For simplicity, let us denote the cohomology class  of a
$D_\mu$-cocycle $\ph$ by $\bar\ph$. It is easily checked that
$D_\mu$ and $D_{\dl}$ anticommute, since $[\mu,\delta]=0$. It
follows that $$D_\mu:\ker(D_{\delta})\ra\ker(D_{\delta}),$$ which
means that we can define the cohomology $H_\mu(\ker(\delta))$,
determined by the restriction of $D_\mu$ to $\ker(D_{\delta})$. The
existence of a $\psi$ such that $[\delta,\lambda]+[\mu,\psi]=0$ and
$[\delta,\psi]=0$ is equivalent to the assertion that
$\overline{[\delta,\lambda]}=0$ in $H_\mu(\ker(\delta))$.

Note that even though the condition for the existence of a $\psi$
depends explicitly on $\lambda$, if the statement is true for
$\lambda$,  then it is  valid for any element in $\bar\lambda$. This
follows because if  $\lambda$ is replaced by
$\lambda'=\lambda-[\mu,\beta]$ and $\psi$ by
$\psi'=\psi-[\delta+\lambda,\beta]$, where $\beta:W\ra M$, then we
obtain a new codifferential $d'=\delta+\mu+t(\lambda'+\psi')$, which
is in fact \emph{infinitesimally equivalent} to $d$. By
infinitesimal equivalence, we mean an equivalence determined  by an
\emph{infinitesimal automorphism} $f=\exp(t\beta)$, where
$\beta:W\ra M$. (Actually, this is a restricted version of
infinitesimal equivalence. We will introduce a more general notion
later.) Since $d'=f^*(d)$, it follows that
$[\delta,\lambda']+[\mu,\psi']=0$.

Now consider a fixed $\lambda$ such that $\overline{[\delta,\lambda]}=0$ in
$H_\mu(\ker(\delta))$, and choose some $\psi$ such  that
$[\delta,\lambda]+[\mu,\psi]=0$.  If $\psi'=\psi+\tau$ is another solution, then $[\mu,\tau]=0$
and $[\delta,\tau]=0$. Now $[\mu,\delta]=0$,so $\bar\delta$ is defined in $\hmu$. The cohomology
$\hmu$ inherits the structure of a graded Lie algebra, with
$$[\bar\alpha,\bar\beta]=\overline{[\alpha,\beta]}.$$
Since $[\bar\delta,\bar\delta]=\overline{[\delta,\delta]}=0$,
$\bar\delta$ determines a coboundary operator
$D_{\bar\delta}$ on $\hmu$. Denote the cohomology of $D_{\bar\delta}$ by $H_{\mu,\delta}$, and the
cohomology class of a $D_{\bar\delta}$-cocycle $\bar\ph$ by $[\bar\ph]$. Then $[\bar\delta,\bar\tau]=0$,
so $\tau$ determines a cohomology class $[\bar\tau]$. On the other hand, suppose that $\bar\tau$ is
any $D_{\bar\delta}$-cocycle.
Then $[\bar\delta,\bar\tau]=0$ implies that $[\delta,\tau]$ is a $D_\mu$-coboundary. Since
$[\delta,\tau]\in C^{0,3}$, this forces $[\delta,\tau]=0$. Thus, every $D_{\bar\delta}$-cocycle determines
an extension.

We need to determine when an extension given by a $D_{\bar\delta}$-cocycle $\bar\tau$ is equivalent to the extension
given by  $\lambda$ and $\psi$. Suppose that $\bar\tau=-D_{\bar\delta}(\bar\beta)$. This is equivalent to
the conditions $[\mu,\beta]=0$ and $\tau=-[\delta,\beta]$,  which is equivalent to the
condition $f^*(d)=d'$, where
$f=\exp(t\beta)$. Thus, infinitesimal equivalence classes  are parameterized by $H_{\mu,\delta}$.
We summarize these results in the following theorem.
\begin{thm}\label{th3}
The infinitesimal extensions of the Lie algebra structure $\delta$ on $W$ by a Lie algebra structure
$\mu$ on $M$ are completely classified by the set of $D_\mu$-cohomology classes $\bar\lambda$
arising from $\lambda\in\hom(MW,M)$ which satisfy the formula
$$\overline{[\delta,\lambda]}=0\in\hrest$$
together with the $D_{\bar\delta}$-cohomology classes $[\bar\tau]$ of $D_{\bar\delta}$-cocycles
$\bar\tau$  arising from of $D_\mu$-cocycles
$\tau\in\hom(W^2,M)$.
\end{thm}

\subsection{Classification of extensions of a Lie algebra by a module}

We consider the special case when the structure $\lambda$ determines an actual module structure on
$W$; in other words, when the MC formula
\begin{equation*}
[\delta,\lambda]+\tfrac12[\lambda,\lambda]=0
\end{equation*}
holds.
In this case it follows that $[\mu,\psi]=0$, so that $\bar\psi$ is well defined.
The condition $[\delta+\lambda,\psi]=0$
is equivalent to the condition $D_{\bdl}(\bar\psi)=0$
where  $D_{\bdl}(\bar\ph)=[\bdl,\bar\ph]$. Note that $D_{\bdl}^2=0$, so if we define the cohomology
class of a $D_{\bdl}$-cocycle $\bar\ph$ by $[\bar\ph]$, then we can express this condition as
$[\bar\psi]=0$.  Moreover, if $D_{\bdl}(\bar\psi)=0$, then $[\delta+\lambda,\psi]=[\mu,\beta]$ for some
$\beta$, but since $[\delta+\lambda,\psi]\in C^{0,3}$, which intersects the
$D_\mu$-cobounaries trivially, this forces $[\delta+\lambda,\psi]=0$. Thus
extensions are given by $D_{\bdl}$-cocycles.

On the other hand, if $\bar\psi=D_{\bdl}(\bar\beta)$, then $\psi=[\delta+\lambda,\beta]$
and $[\mu,\beta]=0$. It follows that the extension determined by $\psi$ is equivalent to the
one given by $\psi=0$.  Thus the equivalence classes of extensions preserving the
module structure $\lambda$ are parameterized  by
$[\bar\psi]$ for $\psi\in\hom(W^2,M)$.

\begin{thm}\label{th4}
The extensions of $\delta$ by $\mu$ determined by a fixed $\lambda$
satisfying
\begin{equation}
[\delta,\lambda]+\tfrac12[\lambda,\lambda]=0.
\end{equation}
are classified by the $D_{\bdl}$-cohomology classes $[\bar\psi]$ of $D_{\bdl}$-cocycles
$\bar\psi$  arising from $D_\mu$-cocycles $\psi\in\hom(W^2,M)$.
\end{thm}
When $\mu=0$, the theorem above is a reformulation, in the language of codifferentials, of the classical theorem
that the equivalence classes of extensions of a Lie algebra $W$ by a module $M$ are
parameterized by $H^2(W,M)$ (see \cite{feigin-fuchs}).
In fact, the cohomology of a Lie algebra with coefficients in a module is determined by the coboundary
operator $D_{\delta+\lambda}$, which coincides with $D_{\bdl}$ when $\mu=0$.
\subsection{Restricted equivalence classes of extensions}
In any extension, we note that $\delta$ and $\lambda$ are $D_\mu$-cocycles. If $\lambda'$ and $\lambda$
are equivalent $D_\mu$-cocycles and the pair $\lambda,\psi$ determines an extension, then there is an equivalent
extension with $\lambda',\psi'$ given by the formulas in \eqref{equivext}.  Moreover, formula
\eqref{mc-twist} yields the MC formula
$$[\bar\delta,\bar\lambda]+\tfrac12[\bar\lambda,\bar\lambda]=0,$$
which means that given an representative $\lambda$ of a cohomology class $\bar\lambda$,
there is a $\psi$ satisfying \eqref{mc-twist} precisely when
$\bar\lambda$ satisfies the MC-equation for
$\bar\delta$, which is a codifferential in $\hmu$.

We also need $\psi$ to satisfy \eqref{ker-cond},
which is not automatic.  The fact that $D_{\dl}$ anticommutes with $D_\mu$ implies that
$D_\mu$ induces a coboundary operator on $\ker(D_{\dl})$. Note that $[\delta+\lambda,\delta+\lambda]$
lies in this kernel, because the triple bracket of any coderivation vanishes. As a consequence, we obtain
that \eqref{mc-twist} is equivalent to the condition that $[\delta+\lambda,\delta+\lambda]$ is a
$D_\mu$-coboundary in the restricted complex $\ker(D_{\dl})$. Moreover, even though the complex
$\ker(D_{\dl})$ depends on $\lambda$,
the existence of a solution to \eqref{mc-twist} depends only on the $D_\mu$-cohomology
class of $\lambda$.  Thus the assertion that $\overline{[\delta+\lambda,\delta+\lambda]}=0$ in
$\hrest$ depends only on $\bar\lambda$, and not on the choice of a representative. However, the $\psi$
satisfying equation \eqref{mc-twist} does depend on $\lambda$. We encountered a similar situation when
analyzing infinitesimal extensions, except that there, one had to consider only $H_\mu(\ker(\delta))$,
instead of $\hrest$.

Now fix $\lambda$ and some $\psi$ satisfying \eqref{mc-twist}.  We want to characterize the set of all
solutions $\psi'=\psi+\tau$, which satisfy \eqref{mc-twist}.  Note that $[\mu,\tau]=0$  and
$[\delta+\lambda,\tau]=0$. This second condition implies that $\bar\tau$ is a $D_{\bdl}$-cocycle.
Moreover, note that $\tau\in C^{0,2}$, so that $[\bdl,\bar\tau]=0$ is equivalent  to $[\mu,\tau]=0$ and
$[\delta+\lambda,\tau]=0$.  On the other hand, if
$\bar\tau=[\bdl,\bar\beta]$, then $[\delta+\lambda,\beta]=0$ and $[\mu,\beta]=0$, and it is immediate
that the extensions given by $\lambda,\psi$ and $\lambda,\psi+\tau$ are equivalent.  Thus $D_{\bdl}$-coboundaries
correspond to trivial extensions. As a consequence, the equivalence classes of nonequivalent extensions
are parameterized by $\hdlm$. Let us denote the cohomology class in $\hdlm$ of a $D_{\bdl}$-cocycle $\bar\ph$
by $[\bar\ph]$. Note that $\hdlm$ is naturally a graded Lie algebra.
We have shown  the following theorem.
\begin{thm}\label{th5}
The equivalence classes of extensions of the Lie algebra structure
$\delta$ on $W$ by a Lie algebra structure $\mu$ on $M$ under the
action of the restricted group of automorphisms $\grest$ are
completely classified by the set of $D_\mu$-cohomology classes
$\bar\lambda$ arising from $\lambda\in\hom(MW,M)$ which satisfy the
MC-equation
$$\overline{[\delta+\lambda,\delta+\lambda]}=0\in\hrest$$
together with the $D_{\bdl}$-cohomology classes $[\bar\tau]$ of
$D_{\bdl}$-cocycles $\bar\tau$ of $D_\mu$-cocycles
$\tau\in\hom(W^2,M)$.
\end{thm}

\subsection{General Equivalence Classes of Extensions}
In the standard construction of equivalence of extensions, we have
assumed that the homomorphism $f:V\ra V$ acts as the identity on $M$
and $W$.  We could consider a more general commutative diagram of
the form
$$
\begin{CD}
0@>>>M@>>>V@>>>W@>>>0\\
@.@VV\eta V@VVfV@VV\gamma V\\
0@>>>M@>>>V@>>>W@>>>0\\
\end{CD}
$$
where $\eta$ and $\gamma$ are isomorphisms.  It is easy to see that
under this circumstance, if $d'$ is the codifferential on the top
line, and $d$ is the one below, then $\eta^*(\mu)=\mu'$ and
$\gamma^*(\delta)=\delta'$.  Therefore, if one is interested in
studying the most general moduli space of all possible extensions of
all codifferentials on $M$ and $W$, where equivalence of elements is
given by diagrams above, then for two extensions to be equivalent,
$\mu'$ must be equivalent to $\mu$ as a codifferential on $M$, and
$\delta'$ must be equivalent to $\delta$ as a codifferential on $W$,
with respect to the action of the automorphism group $\GL(M)$ on $M$
and $\GL(W)$ on $W$.

Thus, in classifying the elements of the moduli space, we first have
to consider equivalence classes of codifferentials on $M$ and $W$.
As a consequence, after making such a choice,  we need only consider
diagrams which preserve $\mu$ and $\delta$;  in other words, we can
assume  that $\eta^*(\mu)=\mu$ and  that $\gamma^*(\delta)=\delta$.

Next note that we can always decompose a general extension diagram
into one of the form
$$
\begin{CD}
0@>>>M@>>>V@>>>W@>>>0\\
@.@|@VVfV@|\\
0@>>>M@>>>V@>>>W@>>>0\\
@.@VV\eta V@VVg=(\eta,\gamma)V@VV\gamma V\\
0@>>>M@>>>V@>>>W@>>>0\\
\end{CD}
$$
where $f=\exp(\beta)$, and $g=(\eta,\gamma)$ is an element of the
group $\gdiag$ consisting of block diagonal matrices. The group
$\ggen$ of general equivalences is just the group of block upper
triangular matrices, and is the semidirect product of $\grest$ with
$\gdiag$; that is, $\ggen=\grest\rtimes\gdiag$. In fact, if
$g\in\gdiag$, then $g\inv\exp(\beta)g=\exp(g^*(\beta))$.

The group $\gdiag$ acts in a simple manner on cochains. If
$g\in\gdiag$, then $g^*(C^{k,l})\subseteq C^{k,l}$ and
$g^*(C^k)\subseteq C^k$. Since $g^*D_\mu=D_{g^*(\mu)}g^*$, the
action induces a map $$g^*:H_\mu\ra H_{g^*(\mu)},$$ given by
$g^*(\bar\ph)=\overline{g^*(\ph)}$. Similarly,
$g^*D_{\dl}=D_{g^*(\delta)+g^*(\lambda)}g^*$, so we obtain a map
$$g^*:H_{\mu,\dl}\ra H_{g^*(\mu),g^*(\delta)+g^*(\lambda)},$$
given by $g([\bar\ph])=[\overline{g^*(\ph)}]$.

Let $\gdiagmd$ be the subgroup of $\gdiag$ consisting of those
elements $g$ satisfying $g^*(\mu)=\mu$ and $g^*(\delta)=\delta$.
Then $\gdiagmd$ acts on $H_\mu$, and induces a map $H_{\mu,\dl}\ra
H_{\mu,\delta+g^*(\lambda)}$. Let $\gdiagmdl$ be the subgroup of $g$
in $\gdiagmd$ such that $g^*(\lambda)=\lambda$. Thus $\gdiagmdl$
acts on both on $H_\mu$ and $H_{\mu,\dl}$.

It is easy to study the behaviour of elements in $\gdiagmd$ on
extensions. If $\lambda$ gives an extension and $g\in\gdiagmd$, then
any element $\lambda'\in g^*(\bar\lambda)$ will determine an
equivalent extension, and thus equivalence classes of $\bar\lambda$
under the action of the group $\gdiagmd$ correspond to equivalent
extensions.

Now suppose that $\lambda$, $\psi$ gives an extension, and
$\bar\tau$ is a $D_{\bdl}$-cocycle. If $g\in\gdiagmdl$, then
$$g^*(\psi+\tau)=\psi+g^*(\psi)-\psi + g^*(\tau),$$ so that
$$[\bar\tau]\mapsto [\overline{g^*(\psi)-\psi +g^*(\tau)}]$$
determines an action of $\gdiagmdl$ on $H_{\mu,\dl}$ whose
equivalence classes determine equivalent representations.

To understand the action of $\ggen$ on extensions, first note that
any element $h\in\ggen$ can be expressed uniquely in the form
$h=g\exp(\beta)$ where $g\in\gdiag$. If $d'=h^*(d)$, for an
extension $d$, then we compute the components of the extension $d'$
as follows.
\begin{align*}
\delta'&=g^*(\delta)\\
\mu'&=g^*(\mu)\\
\lambda'&=g^*(\lambda)+[\mu',\beta]\\
\psi'&=g^*(\psi)+[\delta'+\lambda'-\tfrac12[\mu',\beta],\beta].
\end{align*}
Clearly, $\delta'=\delta$ and $\mu'=\mu$ precisely when
$g\in\gdiagmd$.   Define the group $\ggenmd$ to be the subgroup of
$\ggen$ consisting of those $h=g\exp(\beta)$ such that
$g^*(\delta)=\delta$ and $g^*(\mu)=\mu$. We have a simple
decomposition $\ggenmd=\gdiagmd\rtimes\grest$.

Define $\ggenmdl$ to be the subgroup of $\ggenmd$ consisting of
those $h$ such that $\lambda=g^*(\lambda)+[\mu,\beta]$, then
$\ggenmdl$ does not have a a simple decomposition in terms of
$\gdiagmdl$, because the condition $\lambda'=\lambda$ does not force
$g\in\gdiagmdl$. However, we can still define an action of
$\ggenmdl$ on $H_{\mu,\dl}^{0,2}$ by
$$[\bar\tau]\ra [\overline{
g^*(\psi)-\psi+g^*(\tau)+[\delta+\lambda-\tfrac12[\mu,\beta],\beta]}],$$
whose equivalence classes determine equivalent representations. Note
that for any element $\ph$ in $C^{0,2}$, $g^*(\ph)=h^*(\ph)$, so we
can use $h$ in place of $g$ in the formula above.

\begin{thm}\label{th6}
The equivalence classes of extensions of $W$ by $M$ under the action
of the group $\ggen$ are classified by the following data:
\begin{enumerate}
\item Equivalence classes of codifferentials $\delta$ on $W$ under the action
$\GL(W)$.
\item Equivalence classes of codifferentials
$\mu$ on $M$ under the action of the group $\GL(M)$.
\item Equivalence classes of $D_\mu$-cohomology classes $\bar\lambda\in H_\mu^{1,1}$
which satisfy the MC-equation
$$\overline{[\delta+\lambda,\delta+\lambda]}=0\in H^{1,2}_\mu(\ker(D_{\dl}))$$
under the action of the group $\gdiagmd$ on $\hmu$.
\item Equivalence classes of $D_{\bdl}$-cohomology classes $[\bar\tau]\in H_{\mu,\dl}^{0,2}$
under the action of the group $\ggenmdl$.
\end{enumerate}
\end{thm}
We are more interested in the moduli space of extensions of $W$ by
$M$ preserving fixed codifferentials on these spaces.
\begin{thm}\label{th7}
The equivalence classes of extensions of a codifferential $\delta$
on $W$ by a codifferential $\mu$ on $M$ under the action of the
group $\ggenmd$ are classified by the following data:
\begin{enumerate}
\item Equivalence classes of $D_\mu$-cohomology classes $\bar\lambda\in H_\mu^{1,1}$
which satisfy the MC-equation
$$\overline{[\delta+\lambda,\delta+\lambda]}=0\in H^{1,2}_\mu(\ker(D_{\dl}))$$
under the action of the group $\gdiagmd$ on $\hmu$.
\item Equivalence classes of $D_{\bdl}$-cohomology classes $[\bar\tau]\in H_{\mu,\dl}^{0,2}$
under the action of the group $\ggenmdl$.
\end{enumerate}
\end{thm}
To illustrate why this more general notion of equivalence is useful, we give some simple examples
of extensions of ordinary and super Lie algebras.
\subsection{Simple examples of (super) Lie algebra extensions}
We first study the simplest possible cases of extensions, where $W$ and $M$ are both 1-dimensional.
However, we also have to take into account whether the basis elements are even or odd, which means
that we obtain 3 different nontrivial cases, according to whether both $W$ and $M$ have odd bases,
corresponding to the ordinary Lie algebra situation, or whether one of these spaces has an odd basis.
(If both spaces  have even bases, then there are no nontrivial codifferentials.)
The only possible codifferentials on these spaces are the trivial ones
$\delta=0$ and $\mu=0$.
\subsubsection{Extending a $0|1$-dimensional space by a $0|1$-dimensional space}
Recall that there is only 1 nontrivial Lie algebra structure on a $0|2$-dimensional space
$V=\langle v_1,v_2\rangle$, given by the codifferential $\psa{11}2$. Let $M=\langle v_1\rangle$
and $W=\langle v_2\rangle$.  Also, since $W^2=0$, the cocycle $\psi$ must vanish. The coderivation
$\lambda=a\psa{12}1$ determines an extension for every value of $a$. If one considers only
restricted equivalences, then the extensions for different values of $a$ are not equivalent,
but if one considers the general notion of equivalence, then the codifferentials associated to
nonzero values of $a$ are equivalent.  Note that the moduli space of equivalence classes of
codifferentials on $V$ consists of only one point, corresponding to the codifferential
$\psa{12}1$. Thus the general notion of equivalence gives a more natural connection with this
moduli space. Note that not every codifferential on $V$ arises from an extension of this type,
but rather the equivalence classes correspond.
\subsubsection{Extensions of a $0|1$-dimensional by a $1|0$-dimensional space}
Let $M=\langle v_1\rangle$
and $W=\langle v_2\rangle$. Then any nontrivial extension is equivalent to the extension given
by  $\lambda=\psa{12}1$, and the cocycle $\psi$ must vanish.
\subsubsection{Extensions of a $1|0$-dimensional by a $0|1$-dimensional space}
Let $M=\langle v_2\rangle$
and $W=\langle v_1\rangle$. (Note that we change the order of the basis to conform with the principle of
listing the even basis element of the space $V=\langle  v_1,v_2\rangle$  first.)
This time, the module structure $\lambda$ must vanish, but we do have a nontrivial cocycle $\psi$
which is always equivalent to $\psi=\psa{11}2$.

Let us compare the two moduli spaces we have just studied with the moduli space of $1|1$-dimensional
super Lie algebras. If $d$ is an odd element of $\hom(V^2,V)$, then $d=a\psa{12}1+b\psa{11}2$.
However, $[d,d]=0$ precisely when $a=0$ or $b=0$. Thus, the moduli space consists of exactly two
elements, corresponding to the union of the two moduli spaces of extensions.
\subsubsection{The moduli space of 3-dimensional  Lie algebras}
Let $V$ be a 3-dimensional space with a completely odd basis, so the moduli space
of codifferentials on this space corresponds to the moduli space of ordinary 3-dimensional Lie algebras.
We recall the decomposition of this moduli space into equivalence classes represented by the codifferentials
below as given in \cite{fp4}.
\begin{align*}
d_1&=\psa{23}1\\
d_2&=\psa{13}1+\psa{23}2\\
d_3&=\psa{12}3+\psa{13}2+\psa{23}1\\
d(\lambda:\mu)&=\psa{13}1\lambda+\psa{23}2\mu+\psa{23}1
\end{align*}
The codifferential $d_3$ gives  the simple Lie algebra $\mathfrak{sl}_2(\C)$, which is not given by an extension.
The other codifferentials all arise from extensions of the trivial 1-dimensional Lie algebra by a 2-dimensional
Lie algebra, as well as by  extensions of a 2-dimensional Lie algebra by the trivial Lie algebra. We analyze
the resulting moduli spaces.
\subsubsection{Extensions of a $0|2$-dimensional by a $0|1$-dimensional space}\label{example4}
Let $W=\langle v_1,v_2\rangle$ and $M=\langle v_3\rangle$.
On $M$, there is only one codifferential, $\mu=0$. There are two nonequivalent codifferentials on
$W$, given  by $\delta=\psa{12}1$ and the trivial codifferential $\delta=0$.
The coderivation $\lambda$ must be of the form $\lambda=\psa{13}3a+\psa{23}3b$, and  the cocycle $\psi$
must be of the form $\psi=\psa{12}3c$. For any values of $a$ and $b$, we have $[\lambda,\lambda]=0$,
and $[\delta+\lambda,\psi]=0$, for both choices of $\delta$ and for any values of the parameters $a,b,c$.
 Note that independently
of the choice of $\delta$, $\lambda$ gives a true module structure for any extension, because $\mu=0$.
Therefore $\psi=0$ always gives an extension for any solution $\lambda$ to the MC-equations.

If $\delta=\psa{12}1$, then $[\delta,\lambda]=\psa{123}3a$, and therefore $a=0$.
If $\beta=\pha13x+\pha23y$, then $[\delta+\lambda,\beta]=\psa{12}3(ay-(1+b)x)$. As a consequence,
the $\tau$ term in the classification can be taken as $\tau=\psa{12}3c$. Moreover, unless $a=0$ and
$b=-1$, we can always find values of $x$ and $y$ so that $\tau=[\delta+\lambda,\beta]$. Thus, except
in this special case, $[\bar\tau]=0$, and so we  can assume that $c=0$.
Therefore, we
can assume that $\tau=0$ unless $b=-1$.
In this case, the group
$G(\mu,\delta)$ of automorphisms preserving $\mu$ and $\delta$ is given by automorphisms $g$ whose
 matrices are of the form
$$
\left[
\begin {array}{ccc}
r&s&0\\\noalign{\medskip}0&1&0\\\noalign{\medskip}0&0&t
\end {array}
\right],
$$
where $r,s,t$ are arbitrary parameters such that $rt\ne0$. It is easily checked that $g^*(\lambda)=\lambda$
for every such $g$, so different values of  $b$ give rise to distinct
$D_\mu$-cohomology classes $\bar\lambda=\lambda$
(since $\mu=0$). Note that the MC formula is satisfied by $\lambda$, so $\lambda$ determines a
$\delta$-module structure on $W$. Since $G(\delta,\mu)$
is the automorphism group of $\lambda$-preserving automorphisms of $V$, we have to consider the action of
such an automorphism on $[\bar\tau]$. We compute that $g^*(\tau)=\tau(r/t)$,
which means that when $b=-1$, when we have to consider the case $c\ne0$, we  can choose $c=1$ to represent
the  isomorphism class of $[\bar\tau]$.

Putting  this all together, we obtain the following codifferentials. When $\tau=0$ we obtain the
codifferential $d=\psa{12}1+\psa{23}3b\sim d(1:-b)$,  unless
$b=-1$, in which case we obtain the codifferential $d_2$. When $b=-1$ and $c=1$ we obtain the codifferential
$d=\psa{12}1-\psa{23}3+\psa{12}3\sim d(1:1)$.
Note that the extensions of $\delta$ are never equivalent to $d_1$.

Consider now the case when
$\delta=0$. In this case, $[\delta,\lambda]=0$, so neither $a$ nor $b$  is forced to vanish. The group
$G(\delta,\mu)$ consists of matrices of the form
$$
\left[
\begin {array}{ccc}
r&s&0\\\noalign{\medskip}p&q&0\\\noalign{\medskip}0&0&t
\end {array}
\right],
$$
such that $t(rq-sp)\ne0$. In this case $g^*(\lambda)=\lambda=\psa{13}3(ar+bp)+\psa{23}3(as+bq)$.
Thus, there are 2 isomorphism classes for $\lambda$, represented by $\lambda=\psa{23}3$ and $\lambda=0$.
Again, let  $\beta=\pha13x+\pha23y$.
For $\lambda=\psa{23}3$,  we have $[\delta+\lambda,\tau]=-\psa{12}3x$, so $\bar\tau$ is a
$D_{\bdl}$-coboundary.  Thus, we do not have to consider the action  of the
group of $\lambda$-preserving automorphisms. Thus we have only one codifferential  to consider
$d=\psa{23}3\sim d(0:1)$.

For $\lambda=0$, we cannot assume that $\tau=0$.
The group of $\lambda$-preserving automorphisms
is just $G(\delta,\mu)$, we have two isomorphism classes for $\tau$,
$\tau=0$ and $\tau=\psa{12}3$. If $\tau=0$, then
$d=0$. Otherwise, $d=\psa{12}3$, which is equivalent to the codifferential $d_1$. This was the one case
which did not show up for $\delta=\psa{12}1$.
\subsubsection{Extensions of $0|1$-dimensional by $0|2$-dimensional spaces}
In \cite{fp3}, the moduli space of 3-dimensional Lie  algebras was constructed by considering
extensions of $\C$ by a 2-dimensional Lie algebra, in other words, by exactly the consideration
we now present, although from a slightly different point of view. The calculations we present here
are given in more detail in \cite{fp3,Ott-Pen}.
Let $V=\langle v_1,v_2\rangle$ and $W=\langle v_3\rangle$.
Since $W$ is 1-dimensional $\delta=\psi=0$.
We must have  $\lambda=\psa{13}1a+\psa{23}1b+\psa{13}2x+\psa{23}2y$ for some  values of the  parameters $a,b,x,y$.
In all cases $[\lambda,\lambda]=0$, so  $\lambda$ determines a true module structure. There  are,
up to isomorphism, two possible cases, $\mu=\psa{12}1$ and $\mu=0$.

Consider the case when $\mu=\psa{12}1$. In this case, the condition $[\mu,\lambda]=0$ forces
$x=y=0$. Moreover, if  $\beta=\pha31s+\pha32t$, then $[\mu,\beta]=\psa{13}1t-\psa{23}1s$, which
means that $\lambda$ is always a coboundary, so we need only consider the case $\lambda=0$.
But then the extension is simply given by the codifferential $\mu$, which is equivalent to the codifferential
$d(0:1)$.
The case $\mu=0$ gives no condition on $\lambda$, so we need only consider the action of the group
$G(\mu,\delta)$ on $\lambda$.   This action produces precisely  the similarity classes (up to multiples)
of the $2\times 2$  matrices  $\left[\begin{smallmatrix}a&b\\x&y\end{smallmatrix}\right]$.  This is exactly
the characterization of the codifferentials $d_1$,  $d_2$, and $d(\lambda:\mu)$,
which are given by the matrices
$\left[\begin{smallmatrix}0&1\\0&0\end{smallmatrix}\right]$,
$\left[\begin{smallmatrix}1&0\\0&1\end{smallmatrix}\right]$, and
$\left[\begin{smallmatrix}\lambda&1\\0&\mu\end{smallmatrix}\right]$, respectively. Thus every codifferential
on $V$ except $d_3$ arises as an extension of $\delta=0$ by $\mu=0$.
\subsubsection{A more complicated extension}
In order to illustrate this construction with a more interesting example,
we would like to consider an case where neither $\delta$
nor $\lambda$ vanish, which is not possible if the total space has dimension 3.
We give an example
on a space of total dimension 5.

Let $M=\langle v_1,v_2,v_3\rangle$ and $W=\langle v_4,v_5\rangle$. Let $\mu=d(0:1)=\psa{23}1+\psa{23}2$
and $\delta=\psa{45}4$. If we denote the $5\times 10$ matrix of the extended codifferential
by $A=(a_{ij})$, then $\lambda$ is determined by the submatrix
$$L=
\left[\begin{matrix}
a_{14}&a_{15}&a_{16}&a_{17}&a_{18}&a_{19}\\
a_{24}&a_{25}&a_{26}&a_{27}&a_{28}&a_{29}\\
a_{34}&a_{35}&a_{36}&a_{37}&a_{38}&a_{39}\\
\end{matrix}\right].
$$
Taking into account the condition $[\mu,\lambda]=0$ yields a much simpler matrix for $\lambda$:
$$L=\left[\begin{matrix}
a_{14}&a_{15}&a_{16}&a_{17}&a_{18}&a_{19}\\
0&a_{14}+a_{15}&a_{26}&0&a_{17}+a_{18}&a_{29}\\
0&0&0&0&0&0\\
\end{matrix}\right].
$$
Let $\beta=\pha41c_{11}+\pha42c_{21}+\pha43c_{31}+\pha51c_{12}+\pha52c_{22}+\pha53c_{32}\in\hom(W,M)$.
Then the matrix corresponding to $[\mu,\beta]$ is
$$\left[\begin{matrix}
0&c_{31}&-c_{21}&0&c_{32}&-c_{22}\\
0&c_{31}&-c_{21}&0&c_{32}&-c_{22}\\
0&0&0&0&0&0
\end{matrix}\right],
$$
so by adding a $D_\mu$-coboundary to $\lambda$, we reduce to the case when
$$L=\left[\begin{matrix}
a_{14}&0&a_{16}&a_{17}&0&a_{19}\\
0&a_{14}&0&0&a_{17}&0\\
0&0&0&0&0&0\\
\end{matrix}\right].
$$
If we express $\psi$ in the form $\psi=\psa{45}1b_1+\psa{45}2b_2+\psa{45}3b_3$, then
the condition
$[\delta,\lambda]+\tfrac12[\lambda,\lambda]+[\mu,\psi]=0$ gives $b_2=b_3=0$, $a_{14}=0$ and
either $a_{17}=1$ or $a_{16}=0$.  Note that after taking these conditions into account, we
have $[\mu,\psi]=0$, so that actually, $\lambda$ determines a true module structure on $M$.
Moreover, this means that we can choose $\psi=0$ in satisfying the MC-equation, and that therefore
the cochain $\tau$ in the theorem is just $\tau=\psa{45}1b_1$. It is easily checked that $[\delta+\lambda,\tau]=0$.
Moreover $$[\delta+\lambda,\beta]=\psa{45}1(c_{31}a_{19}-c_{32}a_{16})-\psa{45}3c_{31},$$ which means
that $\bar\tau$ is a $D_{\bdl}$-coboundary unless $a_{16}=0$. In other words,
unless $a_{16}=0$, the extension is a semi-direct product.

Let us now assume that $a_{16}=0$. Thus $\lambda=\psa{15}1+\psa{25}2=\psa{35}1a_{19}$ and
$\tau=\psa{45}1b_1$. We next have to take into account the action of $G(\mu,\delta)$ on the
equivalence classes of $\bar\lambda$. It is easy to show that an element $g$ in this group has
matrix given by
$$
\left[\begin{matrix}
g_{11}&g_{21}-g_{11}&g_{13}&0&0\\
0&g_{22}&g_{23}&0&0\\
0&0&1&0&0\\
0&0&0&g_{44}&g_{45}\\
0&0&0&0&1
\end{matrix}\right].
$$
One calculates that the matrix of $g^*(\lambda)$ is
$$\left[ \begin {array}{cccccccccc}
0&0&0&1&0&{\frac {g_{{1,3}}g_{{2,2}}-g_{{2,3}}g_{{2,2}}+g_{{2,3}}g_{{1,1}}+a_{{1,9}}g_{{2,2}}}{g_{{
2,2}}g_{{1,1}}}}\\\noalign{\medskip}0&0&0&0&1&{\frac {g_{{2,3}
}}{g_{{2,2}}}}\\\noalign{\medskip}0&0&0&0&0&0
\end {array} \right]
$$
Since the matrix of $\lambda$ is
$$L=\left[\begin{matrix}
0&0&0&1&0&a_{19}\\
0&0&0&0&1&0\\
0&0&0&0&0&0\\
\end{matrix}\right],
$$
it follows that $g^*(\lambda)$ lies in the same $D_\mu$-cohomology class as $\lambda$ precisely when
$g_{23}=a_{19}(1-g_{11})+g_{13}$.  However, $g^*(\tau)=\tau\frac{g_{44}}{g_{11}}$, which means that
if $b_1\ne0$, we can  transform it to $b_1=1$.  This means that there is an extension for which
the cohomology class $[\bar\tau]$ cannot be transformed to the zero cohomology class, which is therefore
an extension of $\delta$ by $\mu$  which is not given by a semi-direct product.
\section{Infinitesimal deformations of extensions of Lie algebras}\label{defext}
A natural question that arises when studying the moduli spaces arising from extensions is how to fit the
moduli together as a space, and to answer that question, one needs to have a notion of how to move around
in the moduli space.  This notion is precisely the idea of deformations, in this case, deformations of
the extensions.  We will classify the infinitesimal deformations of an extension.

Suppose that $d=\delta+\mu+\lambda+\psi$ is an extension, and we consider the infinitesimal deformation
of this extension
$$
d_t=d+t(\eta+\zeta),
$$
where $\eta\in\hom(MW,M)$ represents a deformation of the $\lambda$ structure, and $\zeta\in\hom(W^2,M)$
gives a deformation
of the $\psi$ structure. Here, we don't consider deformations which involve deforming the $\delta$ or $\mu$ structure.
The infinitesimal condition  is that $t^2=0$, in which case, as usual, the condition for $\eta,\zeta$ to determine
a deformation is, infinitesimally, that $[d,\eta+\zeta]=0$.
We split this one condition up into the four conditions below.
\begin{align}
&[\delta+\lambda,\eta]+[\mu,\zeta]=0\label{cond1}\\
&[\delta+\lambda,\zeta]+[\psi,\eta]=0\label{cond2}\\
&[\mu,\eta]=0\label{cond3}\\
&[\psi,\zeta]=0\label{cond4}
\end{align}
These conditions are symmetric in the roles of $\psi$ and $\mu$, but as in the case of infinitesimal extensions,
this symmetry is a bit misleading. For example the condition \eqref{cond4} is automatic for
$\zeta\in\hom(W^2,M)$, but condition \eqref{cond3} is not automatic for $\eta\in\hom(MW,M)$.

Note that since $\eta$ is a $D_\mu$-cocycle, $\bar\eta$ is well defined, and the first of these equations implies
that $\bar\eta$ is a $D_{\bdl}$-cocycle.  Since $[\psi,\psi]=0$, it determines a coboundary
operator $D_\psi$ as well.  Denote the $D_\psi$-cohomology class of a $D_\psi$-cocycle $\ph$ by
$\dbar\ph$ and the set of cohomology classes by
$\hpsi$. Note that $\hpsi$ inherits the structure of a Lie algebra.

Since $[\delta+\lambda,\psi]=0$, it follows that $\bbdl$ is well defined.
Moreover, we have
\begin{align*}
[\dbar{\delta+\lambda},[\dbar{\delta+\lambda},\dbar\ph]]&=\dbar{[\delta+\lambda,[\delta+\lambda,\ph]]}=
\dbar{[\tfrac12[\delta+\lambda,\delta+\lambda],\ph]}
\\&=-\dbar{[[\mu,\psi],\ph]}=-\dbar{[\mu,[\psi,\ph]]}=0,
\end{align*}
so $D_{\bbdl}$ is a differential on $\hpsi$.
Denote the cohomology class of a $D_{\bbdl}$-cocycle $\dbar\ph$ by
$\cbdlp\ph$ and the set of cohomology classes by $\hdlp$. Note that $\hdlp$ inherits the structure
of a Lie algebra.

We first remark that conditions \eqref{cond1} and \eqref{cond3} imply that
$\left[\bar\eta\right]$ is well defined, and \eqref{cond2} and \eqref{cond4} imply that
$\cbdlp\zeta$ is well defined.

Next we introduce an action of $D_\psi$ on $\hdlm$.
It is not possible to extend the operation of bracketing with $\psi$ to the $D_\mu$-cohomology,
because $[\mu,\psi]\ne0$.  Moreover, even if $[\mu,\ph]=0$, it does not follow that $[\mu,[\psi,\ph]]=0$.
However, we can  extend  the bracket to $\hdlm$ as follows.
A cohomology class $\left[\bar\ph\right]$, is given by a $\ph$ such that $[\mu,\ph]=0$ and
$[\delta+\lambda,\ph]=[\mu,\beta]$ for some $\beta$. Note that
\begin{align*}
[\mu,[\psi,\ph]]&=[[\mu,\psi],\ph]=-[\delta+\lambda,[\delta+\lambda,\ph]]
\\&=-[\delta+\lambda,[\mu,\beta]]=[\mu,[\delta+\lambda,\beta]].
\end{align*}
This suggests that it might be possible to define a map  by
$$D_\psi([\bar\ph])=\left[\overline{[\psi,\ph]-[\delta+\lambda,\beta]}\right],$$
where $\beta$ is any solution to $[\delta+\lambda,\ph]=[\mu,\beta]$.
There are a few things we have to check in order to  see that this is a well-defined
action.
\begin{enumerate}
\item Show that $[\bdl,\overline{[\psi,\ph]-[\delta+\lambda,\beta]}]=0$. This follows
from
\begin{align*}
[\delta+\lambda,[\psi,\ph]-[\delta+\lambda,\beta]]&=
[\delta+\lambda,[\psi,\ph]]-[\delta+\lambda,[\delta+\lambda,\beta]]
\\&=-[\psi,[\delta+\lambda,\ph]]+[[\mu,\psi],\beta]
\\&=-[\psi,[\mu,\beta]]+[[\mu,\psi],\beta]=
[\mu,[\psi,\beta]].
\end{align*}
This first step has established that $\left[\overline{[\psi,\ph]-[\delta+\lambda,\beta]}\right]$ is at
least well-defined.
\item Show that, given a $\ph$, the formula does not depend on the choice of $\beta$.  Suppose that
$[\delta+\lambda,\ph]=[\mu,\beta']$. Then $[\mu,\beta'-\beta]=0$
\begin{align*}
&\left[\overline{[\psi,\ph]-[\delta+\lambda,\beta]}-\overline{[\psi,\ph]-[\delta+\lambda,\beta']}\right]\\=
&\left[\overline{[\delta+\lambda,\beta'-\beta]}\right]=\left[[\bdl,\overline{\beta'-\beta}]\right]=0,
\end{align*}
because $\bdl$ is a differential on $\hmu$.
\item Show that if $[\bar\ph]=0$,  then $\left[\overline{[\psi,\ph]-[\delta+\lambda,\beta]}\right]=0$.
If  $[\bar\ph]=0$,  then we can express  $\ph=[\delta+\lambda,\alpha]+[\mu,\gamma]$ for some $\alpha$ and $\gamma$ such
that $[\mu,\alpha]=0$.
Now
\begin{align*}
[\psi,\ph]-[\delta+\lambda,\beta]&=[\psi,[\delta+\lambda,\alpha]+[\mu,\gamma]]-[\delta+\lambda,\beta]\\&=
-[\delta+\lambda,[\psi,\alpha]+\beta]+[\psi,[\mu,\gamma]\\&=
-[\delta+\lambda,[\psi,\alpha]+\beta]+[[\psi,\mu],\gamma]-[\mu,[\psi,\gamma]]\\&=
-[\delta+\lambda,[\psi,\alpha]+\beta]-[\delta+\lambda,[\delta+\lambda,\gamma]]-[\mu[\psi,\gamma]]\\&=
-[\delta+\lambda,[\psi,\alpha]+\beta+[\delta+\lambda,\gamma]]-[\mu,[\psi,\gamma]].
\end{align*}
The last term above is a $D_{\bdl}$-coboundary, provided  that $$[\mu,[\psi,\alpha]+\beta+[\delta+\lambda,\gamma]]=0.$$
But
\begin{align*}
[\mu,[\psi,\alpha]+\beta+[\delta+\lambda,\gamma]]&=-[\mu,[\psi,\alpha]]+[\mu,\beta]+[\mu,[\delta+\lambda,\gamma]]
\\&\hspace{-.3in}=-[\delta+\lambda,[\delta+\lambda,\alpha]]+[\delta+\lambda,\ph]-[\delta+\lambda,[\mu,\gamma]]
\\&\hspace{-.3in}=[\delta+\lambda,-[\delta+\lambda,\alpha]+\ph-[\mu,\gamma]]=0.
\end{align*}
\end{enumerate}
It is straightforward to verify  that $D_\psi$ is actually a Lie algebra morphism on the graded Lie algebra
$\hdlm$, and that $D_\psi^2=0$. Let us denote the cohomology by $\hdlmp$, and
the $D_\psi$-cohomology class of a $D_\psi$-cocycle  $\cbdl{\ph}$ by $\ctpsi{\ph}$.

In a very similar manner, one can show that one can define $D_\mu$  on
$\hdlp$ by
$$
D_\mu(\cbdlp\ph)=\cbdlp{[\mu,\ph]-[\delta+\lambda,\beta]},
$$
where $\beta$ is any coderivation satisfying $[\delta+\lambda,\ph]=[\psi,\beta]$. Then $D_\mu$ is a Lie
algebra morphism on $\hdlp$ whose square is zero, and we denote the resulting  cohomology by $\hdlpm$
and the cohomology class of a $D_\mu$-cocycle $\cbdlp{\ph}$ by $\ctmu{\ph}$.

Thus we have constructed two different triple cohomology groups. It turns out that the first one will play
a more important role in the classification of infinitesimal deformations of extensions.
\begin{lma}\label{eta-zeta}
Suppose that $d=\delta+\mu+\lambda+\psi$ is an extension of the codifferentials $\delta$ on $W$ by
$\mu$ on $M$, that $\eta\in\hom(MW,M)$ and $\zeta\in\hom(W^2,M)$.
If
$$d_t=d+t(\eta+\zeta)$$
determines an infinitesimal deformation of $d$ then
\begin{enumerate}
\item $\ctpsi\eta$ is well defined.
\item $\ctmu\zeta$ is well defined.
\end{enumerate}
\end{lma}
\begin{proof}
We already saw that if $\eta$ and $\zeta$ determine a deformation, then $\cbdl\eta$ and
$\cbdlp\zeta$ are defined.  Now, $D_\psi(\cbdl\eta)=\cbdl{[\psi,\eta]-[\delta+\lambda,\beta]}$,
where $[\delta+\lambda,\eta]=[\mu,\beta]$. We obtain that
$$
[\psi,\eta]-[\delta+\lambda,\beta]=-[\delta+\lambda,\zeta]-[\delta+\lambda,\beta]=-[\delta+\lambda,\beta+\zeta].
$$
Moreover $$[\mu,\beta+\zeta]=[\mu,\beta]+[\mu,\zeta]=[\delta+\lambda,\eta]-[\delta+\lambda,\eta]=0.$$
Thus $\overline{[\psi,\eta]-[\delta+\lambda,\beta]}=[\bdl,-\overline{\alpha+\beta}]$,
so that $\cbdl{[\psi,\eta]-[\delta+\lambda,\beta]}=0$.
The proof that $D_\mu\left(\cbdlp\zeta\right)=0$ is similar.
\end{proof}
 Now, let us suppose that $\ctpsi{\eta}$ exists.
We will show that this condition alone is sufficient to guarantee the existence of a $\zeta$ such that
conditions \eqref{cond1}, \eqref{cond2}, \eqref{cond3} and \eqref{cond4} are satisfied.
Since $[\bar\eta]$ is defined, there is some $\beta\in\hom(W^2,M)$ such that
\begin{equation*}
[\delta+\lambda,\eta]=[\mu,\beta].
\end{equation*}
Note that this $\beta$ is determined up to a $D_\mu$-cocycle. Now,
\begin{equation*}
0=D_\psi([\bar\eta])=[\overline{[\psi,\eta]-[\delta+\lambda,\beta]}],
\end{equation*}
so
\begin{equation*}
\overline{[\psi,\eta]-[\delta+\lambda,\beta]}=[\bdl,\bar\alpha],
\end{equation*}
and therefore
\begin{equation*}
[\psi,\eta]-[\delta+\lambda,\beta]=[\delta+\lambda,\alpha]+[\mu,\gamma].
\end{equation*}
Since $[\psi,\eta]-[\delta+\lambda,\beta]\in\hom(W^3,M)$, the term $[\mu,\gamma]$ must vanish.
Now, by construction, $[\mu,\alpha]=0$, so we observe that
\begin{align*}
[\delta+\lambda,\eta]&=[\mu,\beta+\alpha]\\
[\psi,\eta]&=[\delta+\lambda,\beta+\alpha]
\end{align*}
But this means that $\zeta=-(\beta+\alpha)$ satisfies conditions \eqref{cond1} and \eqref{cond2}, and thus
all the conditions are satisfied, since condition \eqref{cond3} is that $\bar\eta$ exists
and condition \eqref{cond4} is automatic.
In other words, we have constructed a solution to the infinitesimal deformation
problem, given an arbitrary $D_\psi$-cocycle $\cbdl\eta$.

Next we ask what variation is possible in $\zeta$;
in other words, given that $\zeta$ is a solution, when is
$\zeta'=\zeta+\tau$ another solution.
It is easy to see that $[\mu,\tau]=0$ and $[\delta+\lambda,\tau]=0$.
This means that $\bar\tau$ is a $D_{\bdl}$-cocycle,
so that the cohomology class  $\cbdl\tau$ is defined.
Moreover, $[\psi,\tau]=0$ automatically, so it follows that $\ctpsi\tau$ is defined.
Note that for any $\tau\in\hom(W^2,M)$, $\cbdl\tau$ is defined precisely when the conditions
$[\mu,\tau]=0$ and $[\delta+\lambda,\tau]=0$ are satisfied,
and the existence of $\cbdl\tau$ implies the
existence of $\ctpsi\tau$.

In light of lemma \ref{eta-zeta}, one might expect that $\ctmu\tau$ would be the natural
object of study, and indeed, the conditions on $\tau$ imply that $\ctmu\tau$ exists.
However, the existence of $\ctmu\tau$ is not equivalent to the conditions $[\mu,\tau]=0$
and $[\delta+\lambda,\tau]=0$, so it turns out that $\ctmu\tau$ is not the right object
to classify the infinitesimal deformations.

To classify the infinitesimal deformations, we need to consider infinitesimal equivalences.
These are given by maps $g=\exp(t(\alpha+\beta+\gamma))$, where $\alpha\in\hom(M,M)$,
$\beta\in\hom(W,M)$ and $\gamma\in\hom(W,W)$, and we suppose that $t^2=0$. If we
set $d_t'=g^*(d_t)$, then $d_t'=d+t(\eta'+\zeta')$, where
\begin{align*}
\eta'&=\eta+[\delta+\lambda,\alpha+\gamma]+[\mu,\beta]\\
\zeta'&=\zeta+[\psi,\alpha+\gamma]+[\delta+\lambda,\beta].
\end{align*}
Moreover, $[\delta,\alpha+\gamma]=[\mu,\alpha+\gamma]=0$, because we require that
$g$ preserves $\delta$ and $\mu$.

Let us examine the conditions for two cohomology classes $\ctpsi{\eta}$ and $\ctpsi{\eta'}$
to be equal. First, we must have
\begin{equation*}
\cbdl{\eta'}=\cbdl{\eta}+D_\psi(\cbdl\sigma).
\end{equation*}
However, the only way this could be possible is if $\sigma\in\hom(M,W)$, which is not allowed.
Thus, when $\eta\in\hom(MW,M)$, the cohomology class $\ctpsi\eta$ coincides with the
cohomology class $\cbdl\eta$. However, the existence of $\cbdl\eta$ does not imply
the existence of $\ctpsi\eta$ as it did for $\tau\in\hom(W^2,M)$.
Continuing our analysis, we obtain
\begin{equation*}
\overline{\eta'}=\bar\eta+[\bdl,\bar\xi],
\end{equation*}
for some $\xi=\alpha+\gamma$, with $\alpha\in\hom(M,M)$ and $\gamma\in\hom(W,W)$,
such that $[\mu,\alpha+\gamma]=0$. Finally, we can express
\begin{equation*}
\eta'=\eta+[\delta+\lambda,\alpha+\gamma]+[\mu,\beta],
\end{equation*}
where $\beta\in\hom(W,M)$. But this is precisely the $\eta'$ that arises from $g^*(d_t)$. Thus,
in order for two infinitesimal deformations $d_t$ and $d_t'$ to be equivalent, $\eta$ and
$\eta'$ must belong to the same cohomology class in $\hdlmp$.

Next, we can suppose that $\eta'=\eta$, in other words, that
\begin{equation*}
[\delta+\lambda,\alpha+\gamma]+[\mu,\beta]=0.
\end{equation*}
This is precisely the condition that $\cbdl{\alpha+\gamma}$ exists. One should consider the
action of $g^*$ on $d_t$ as follows. If the term from $\hom(W^2,M)$ in $d_t$
is given by $\zeta+\tau$,
then the corresponding term in $d_t'$ should be $\zeta+\tau'$. In other words, we have
\begin{equation*}
\tau'=\tau+[\psi,\alpha+\gamma]+[\delta+\lambda,\beta].
\end{equation*}
But this means that
\begin{equation*}
\cbdl{\tau'}=\cbdl{\tau}+D_\psi(\cbdl{\alpha+\gamma}).
\end{equation*}
Thus $\ctpsi{\tau'}=\ctpsi{\tau'}$.

Now, let us investigate when it is possible for $\ctpsi\tau=\ctpsi{\tau'}$, for
$\tau,\tau'\in\hom(W^2,M)$. First, we have
\begin{equation*}
\cbdl{\tau'}=\cbdl\tau+D_\psi(\cbdl{\sigma})=\cbdl{\tau+[\psi,\sigma]-[\delta+\lambda,\rho]},
\end{equation*}
where
\begin{equation*}
[\delta+\lambda,\sigma]=[\mu,\rho].
\end{equation*}
We obtain that
\begin{equation*}
\overline{\tau'}=\overline{\tau+[\psi,\sigma]-[\delta+\lambda,\rho]}+[\bdl,\bar\xi],
\end{equation*}
for some $\xi$. Finally, we obtain
\begin{equation*}
\tau'=\tau+[\psi,\sigma]-[\delta+\lambda,\rho+\xi]+[\mu,\epsilon].
\end{equation*}
However, because $\tau'\in\hom(W^2,M)$, we must have $[\mu,\epsilon]$. Also, since
$[\mu,\xi]=0$, we have $[\delta+\lambda,\sigma]=[\mu,\rho+\xi]$, so
the combined term $\rho+\xi$ plays the same role for $\sigma$ as $\rho$ does in terms
of the triple coboundary operator $D_\psi$. If we set $\sigma=\alpha+\gamma$ and
$\rho+\xi=-\beta$, then we see that $\cbdl{\tau'}=\cbdl{\tau}+D_\psi(\cbdl{\alpha+\gamma})$.
In other words, the element $\ctpsi{\tau}$ completes the classification of $d_t$ up
to equivalence.  We have shown that

Summarizing these results we obtain
\begin{thm}\label{th8}
Suppose that $d=\delta+\mu+\lambda+\psi$ is an extension of the codifferentials $\delta$ on $W$ by
$\mu$ on $M$.

An element $\eta\in\hom(MW,M)$ gives rise to an infinitesimal deformation for some $\zeta\in\hom(W^2,M)$
if and only if the triple cohomology class $\{[\bar\eta]\}$ is well defined. In this case, if $\zeta$ is
any coderivation such that $\eta$, $\zeta$ determine an infinitesimal deformation, then
$\zeta'=\zeta+\tau$ determines another infinitesimal deformation if and only if the double cohomology
class $[\bar\tau]$ is well defined.

Moreover the infinitesimal equivalence classes of infinitesimal deformations are classified by
the triple cohomology
classes $\{[\bar\eta]\}$ and $\{[\bar\tau]\}$.
\end{thm}
\subsection{An example of an infinitesimal deformation of an extension}
Consider the extensions from the example in section \ref{example4}
Recall that $\delta=\psa{12}1$, $\mu=0$ and $\lambda=\psa{23}3b$.  Let us first consider the case when
$\psi=0$. We must have $\eta=\psa{13}3r+\psa{23}3s$ for some value of the parameters $r$ and $s$. Since
$\mu=0$, $\bar\eta$ is obviously well defined. Next, we have $[\delta+\lambda,\eta]=\psa{123}3r$, so in
order for $[\bar\eta]$ to exist, $r=0$.  Let
\begin{align*}
\alpha&=\pha11a_{11}+\pha12a_{12}+\pha21a_{21}+\pha22a_{22}\\
\beta&=\pha13x+\pha23y\\
\gamma&=\pha33c.
\end{align*}
Then
\begin{equation*}
[\delta+\lambda,\alpha+\gamma]=\psa{12}1a_{22}-\psa{12}2a_{21}+\psa{13}3a_{21}b+\psa{23}3a_{22}b,
\end{equation*}
from which it follows that $\bar\eta$ is not a $D_{\bdl}$-cohomology class. Since $\psi=0$, $\{]\bar\eta]\}$
is defined and is not a $D_\psi$-coboundary.  Therefore, $\{[\eta]\}\ne0$ unless $s=0$.
Next, we must have $\tau=\psa{12}3z$ for some $z$. But then $[\bdl,\bar\tau]=0$, so $\{[\bar\tau]\}$ is defined.
However $[\delta+\lambda,\beta]=-\psa{12}3x(1+b)$, which means that unless $b=-1$, $\bar\tau$ is a
$D_\psi$-coboundary, and therefore $\{[\bar\tau]\}=0$. As a consequence, when $b\ne-1$, the infinitesimal
deformations are given by the $\eta$'s alone.

Now, when $b=-1$ and $\psi=0$, we obtain a nontrivial cohomology class $\{[\bar\tau]\}$, which means that the
infinitesimal deformations are governed by 2 parameters, $s$ and $z$.  Let us consider the other nontrivial case,
when $b=-1$ and $\psi=\psa{12}3$.  It still turns out that $D_\psi(\bar\tau)=0$, but this time, we compute
$[\psi,\alpha+\gamma]=\psa{12}3(a_{11}+a_{22}-c)$. We also must have $[\delta+\lambda,\alpha+\gamma]=0$,
which forces $a_{21}=a_{22}=0$, but even so, there is a solution for $\tau=[\psi,\alpha+\gamma]$. This means
that $\{[\bar\tau]\}=0$, and like the case when $b\ne-1$, we only get one parameter for our infinitesimal
deformations.

In \cite{fp3,Ott-Pen}, the moduli space of 3-dimensional Lie algebras was studied in detail, and it
was shown that codifferentials of the form $d(1:b)$ deform only along this family, but the codifferential
$d_2$ has a jump deformation to $d(1:1)$ as well as deformations along the family. The picture we have
painted of the infinitesimal deformations of  the extensions  mirrors that behaviour.

In fact, there are natural notions of  universal infinitesimal deformations as well as miniversal deformations
of extensions, but these ideas deserve another paper.
\section{Infinitesimal Deformations of Representations}\label{infdef}
The notion of deformations of representations has applications in physics, but is not
well studied in the mathematics literature. There is a nice paper \cite{Lev-Sen} which
studies some aspects of these issues. In this section,
we give a complete classification of infinitesimal
deformations of representations of Lie algebras and Lie superalgebras.

Suppose that $M$ is a Lie algebra, with multiplication $\mu$ which
is also a module over $W$. In other words,
we are studying an extension of $W$ by $M$ for which the cocycle $\psi$ vanishes.
There are two interesting problems we could study.
\begin{enumerate}
\item Allow the module structure $\lambda$ and the algebra structure $\delta$ to vary,
but keep $\mu$ fixed. This case includes the study of deformations of a module structure
where the module does not have an algebra structure.
\item Allow the module structure $\lambda$ and the multiplication $\mu$ to vary,
but keep the algebra structure $\delta$ fixed.
\end{enumerate}
In both of these scenarios, we think of the structures on $M$ and $W$ as being distinct,
with interaction only through $\lambda$, so when considering automorphisms of the
structures, it is reasonable to restrict to automorphisms of $S(V)$ which do not mix
the $W$ and $M$ terms, in other words, those automorphisms which are given by an
automorphism of $M$ and and automorphism of $W$.

Then we have the following maps:
\begin{align*}
&D_{\delta}:C^n\ra C^{n+1}\\
&D_{\lambda}:C^n\ra C^{1,n}\\
&D_{\delta+\lambda}:C^{k,l}\ra C^{k,l+1}\\
&D_\mu:C^{k,l}\ra C^{k+1,l}.
\end{align*}
In the setup of this problem, we only are interested in $C^{k,l}$ for $k\ge 1$, so
we shall restrict our space of cochains in this manner.  Because of this restriction,
we note that an element in $C^{1,1}$ can be a $D_\mu$-cocycle, but never a $D_{\mu}$-
coboundary. Moreover $C^n\subseteq\ker(D_\mu)$, so an element in $C^2$ is always a
$D_\mu$-cocycle, and never a $D_\mu$-coboundary.

Because $\psi=0$, the MC-equation $[\delta,\lambda]+\tfrac12[\lambda\lambda]=0$ is
satisfied, so that $D_{\delta+\lambda}^2=0$.
Since
\begin{align*}
(D_\lambda D_\delta+D_{\delta+\lambda}D_\lambda)(\ph)=&
[\lambda,[\delta,\ph]]+[\delta+\lambda,[\lambda,\ph]]\\
=&[[\lambda,\delta],\ph]-[\delta,[\lambda,\ph]]+[\delta,[\lambda,\ph]]+[\lambda,[\lambda,\ph]]
\\=&[[\delta,\lambda],\ph]+[\tfrac12[\lambda,\lambda],\ph]=0.
\end{align*}
 we have
\begin{equation*}
D_\lambda D_\delta+D_{\delta+\lambda}D_\lambda=0.
\end{equation*}
If we denote the $D_\mu$-cohomology class of a $D_\mu$-cocycle $\ph$ by $\bar\ph$ as
usual, then since $\bar\lambda$ and $\bar\delta$ are defined, we get the following
version of this equation, applicable to the cohomology space $H_\mu$.\begin{equation*}
D_{\bar\lambda} D_{\bar\delta}+D_{\bdl}D_{\bar\lambda}=0.
\end{equation*}
As usual, let us denote the $D_{\bdl}$-cohomology class of a $D_{\bdl}$-cocycle $\bar\ph$
by $[\bar\ph]$.

Let us study the first scenario, where we allow $\lambda$ and $\delta$ to vary, in other
words, we consider
\begin{equation*}
d_t=d+t(\delta_1+\lambda_1),
\end{equation*}
where $\delta_1\in C^2$ and $\lambda_1\in C^{1,1}$ represent the variations in $\delta$ and
$\lambda$.  The infinitesimal condition $[d_t,d_t]=0$ is equivalent to the three
conditions for a deformation of a module structure:
\begin{align*}
[\delta,\delta_1]&=0\\
[\lambda,\delta_1]+[\delta+\lambda,\lambda_1]&=0\\
[\mu,\lambda_1]&=0.
\end{align*}
By the third condition above $\bar\lambda_1$ is well defined, and $\bar\delta_1$ is defined.
We claim that if $D_{\bar\delta}(\bar\delta_1)=0$, which is the first condition, then
the $D_{\bdl}$-cohomology class $[D_{\bar\lambda}(\delta_1)]$ is well defined and depends
only on the $D_\delta$-cohomology class of $\delta_1$. It is well defined because
\begin{equation*}
D_{\bdl}D_{\bar\lambda}(\bar\delta_1)=-D_{\bar\lambda}D_{\bar\delta}(\bar\delta_1)=0.
\end{equation*}
To see that it depends only on the $D_\delta$-cohomology
class of $\bar\delta$, we apply $D_{\bar\lambda}$ to a $D_{\bar\delta}$-coboundary
$D_{\bar\lambda}(\bar\ph)$ to obtain
\begin{equation*}
D_{\bar\lambda}D_{\bar\delta}(\bar\ph)=D_{\bdl}D_{\bar\lambda}(-\bar\ph),
\end{equation*}
which is a $D_{\bdl}$-coboundary.
The second condition for a deformation of the module structure implies that
$[D_{\bar\lambda}(\bar\delta)]=0$. Moreover, if this statement holds, then there is
some $\lambda_1$ such that $\delta_1$ and $\lambda_1$ determine a deformation of
the module structure.  We see that $\lambda_1'=\lambda+\tau$ is another solution precisely when
$\bar\tau$ exists and $D_{\bdl}(\bar\tau)=0$. Thus, given one solution
$\lambda_1$, the set of solutions is
determined by the $D_{\bdl}$-cocycles $\tau\in C^{1,1}$.

Now let us consider infinitesimal equivalence.  We suppose that $\alpha\in C^{1,0}$ and
$\gamma\in C^1$, and $g=\exp(t(\alpha+\beta)$. If $d_t'=g^*(\d_t)$ is given by
the cochains $\delta_1'$ and $\lambda_1'$, then we have
\begin{align*}
\delta_1'&=\delta_1+D_\delta(\gamma)\\
\lambda_1'&=\lambda_1+D_{\lambda}(\alpha+\gamma)\\
D_{\mu}(\alpha+\gamma)&=0.
\end{align*}
It follows that the set of equivalence classes of deformations are determined by
$D_\delta$ cohomology classes of $\delta_1\in C^2$. If we fix $\delta_1$ such that
$D_\delta(\delta_1)=0$ and $\lambda_1$ satisfying the rest of the conditions of a
deformation, then expressing $\tau'=\tau+D_{\lambda}(\alpha+\gamma)$.
But, since $D_\delta(\alpha+\gamma)=0$, this means
we can express $\tau'=\tau+D_{\delta+\lambda}(\alpha+\gamma)$ and $D_\mu(\alpha+\gamma)$,
which means that $\bar\tau'=\bar\tau+D_{\bdl}(\bar\alpha+\bar\gamma)$, and the solutions
for $\tau$ are given by $D_{\bdl}$-cohomology classes of
$D_\mu$-cocycles $\tau\in C^{1,1}$.
Thus we obtain
\begin{thm}\label{th9}
The infinitesimal deformations of a module $M$ with multiplication
$\mu$ over a Lie algebra
$\delta$, allowing the algebra structure $\delta$ on $W$ and module structure $\lambda$
to vary are classified by
\begin{enumerate}
\item $D_{\delta}$-cohomology classes of $D_\delta$-cocycles
$\delta_1\in C^2$ satisfying the condition
$$[D_{\bar\lambda}(\bar\delta_1)]=0$$.
\item $D_{\bdl}$-cohomology classes $[\bar\tau]$ of $D_{\bdl}$-cocycles $\bar\tau$ of
$D_\mu$-cocycles $\tau\in C^{1,1}$.
\end{enumerate}
\end{thm}

Finally, let us study the second scenario, where we allow $\lambda$ and $\mu$,
but not $\delta$, to vary.  We write $d_t=d+t(\lambda_1+\mu_1)$, where
$\lambda_1\in C^{1,1}$ is the variation of $\lambda$ and $\mu_1\in C^{2,0}$ is
the variation in $\mu$.  The Jacobi identity $[d_t,d_t]=0$ gives three conditions
for a deformation of the module structure.
\begin{align*}
&D_\mu(\mu_1)=0\\
&D_{\delta+\lambda}(\mu_1)+D_\mu(\lambda_1)=0\\
&D_{\delta+\lambda}(\lambda_1)=0.
\end{align*}
Recall that $D_\mu$ maps $\ker(D_{\delta+\lambda})$ to itself, so
$H_\mu(\ker(\delta+\lambda))$ is well defined. The first condition on a deformation
says that $\bar\mu_1$ is well defined.  We claim that in that case,
$\overline{D_{\delta+\lambda}(\mu_1)}$ is
a well defined element of $H_\mu(\ker(\delta+\lambda)$ which depends only on $\bar\mu_1$.
This is clear, since $D_{\delta+\lambda}(\mu_1)\in\ker(D_{\delta+\lambda})$,
and  $D_\mu D_{\delta+\lambda}(\mu_1)=-
D_{\delta+\lambda}D_u(\mu_1)=0$. The second condition on a deformation says simply that
$\overline{D_{\delta+\lambda}(\mu_1)}=0$, and the fact that this statement is true
in $H_\mu(\ker(\delta+\lambda))$ is the third condition.  Therefore, assuming that
$\overline{D_{\delta+\lambda}(\mu_1)}=0$, we can find a $\lambda_1$ so that all
of the conditions for a deformation are satisfied.

If $\lambda_1+\tau$ gives another solution, then $D_\mu(\tau)=0$
and $D_{\delta+\lambda}(\tau)=0$.
Because we do not allow elements of $C^{0,1}$ as cochains, these two equalities are
equivalent to  $[\bar\tau]$ being well defined.

If $d_t'=\exp(t(\alpha+\gamma))$, then we obtain the following.
\begin{align*}
[\delta,\alpha+\gamma]&=0\\
\lambda_1'&=\lambda+[\lambda,\alpha+\gamma]\\
\mu_1'&=\mu_1+[\mu,\alpha+\gamma].
\end{align*}

Thus, up to equivalence a deformation is given by a $D_\mu$-cohomology class $\bar\mu_1$. If
we fix $\mu_1$, and then look at the variation in $\tau$, one also obtains that up to
equivalence, the deformation is determined by the $D_{\bdl}$-cohomology class
$[\bar\tau]$ of the $D_\mu$-cohomology class $\bar\tau$. Thus we have shown
\begin{thm}\label{th10}
The infinitesimal deformations of a module $M$ with Lie algebra structure
$\mu$ over a Lie algebra
$\delta$, allowing the algebra structure $\mu$ on $W$ and module structure $\lambda$
to vary are classified by
\begin{enumerate}
\item $D_\mu$ cohomology classes $\bar\mu_1$ of $D_\mu$-cocycles $\mu_1$ lying in
$C^{2,0}$.
\item $D_{\bdl}$-cohomology classes $[\bar\tau]$ of $D_{\bdl}$-cocycles $\bar\tau$ of $D_\mu$-cocycles
$\tau$ lying in $C^{1,1}$.
\end{enumerate}
\end{thm}
\section{Conclusions}
The authors original intent was to give a brief review of the theory of extensions of Lie algebras
and then extend these notions to the \linf\ algebra case.  However, it turned out that in the course
of translating the relevant ideas about extensions to the language of codifferentials, we discovered
that some of the material is not readily accessible in the literature.  For example,
we could not locate a reference
for \thmref{th5}, although it seems to be fundamental in the classification of extensions.  The language
of codifferentials, with its emphasis on the \zt-graded Lie algebra structure of the cochains, is natural
for the description of extensions, and therefore, the results seem to us to be easier to see in this
formulation.

An extension is determined by four coderivations, and the underlying theme of this paper was to study
what occurs when two of them are fixed and the other two allowed to vary.  In each case, we were able
to give a cohomological condition which one of the two
structures must satisfy in order that an extension, or deformation
of an extension, existed.  Then, given any one of the second type of structures yielding an extension,
the variation in the second was also given by a purely cohomological condition. It was necessary to
cast things in terms of the variation in the second structure because the set of second structures
which give a solution in terms of the first structure form an affine space.

In studying deformations of extensions, we looked at three problems, the deformation of an extension
given by varying the module and cocycle structures, and deformations given by a module structure,
either by holding the algebra structure on the module fixed and varying the module structure and the
algebra structure on the quotient, or by holding the algebra structure on the quotient fixed and varying
the algebra structure on the module, as well as the module structure.  All three of these problems could
have been solved by simply pointing out that an infinitesimal deformation of an extension given by a
codifferential $d$ is classified by the cohomology $H(d)$.  But this approach would not give separate
criteria on the structures we allow to vary, so although correct, it is not the answer we were looking
for.

In \cite{fp3,fp8}, moduli spaces of 3 and 4 dimensional Lie algebras were studied. Our basic approach
to the classification of the Lie algebras was to construct Lie algebras of higher dimension as extensions
of Lie algebras of lower dimension. Of course, this is the fundamental idea behind the classical
decomposition theory of Lie algebras, so our approach is nothing new in this regard.  What we discovered
that was a new perspective was that by decomposing the moduli space in a careful manner, one could show
that the moduli spaces have natural stratifications by orbifolds.  The manner in which we constructed
the moduli spaces by extensions was quite helpful in discovering that stratification.




\providecommand{\bysame}{\leavevmode\hbox to3em{\hrulefill}\thinspace}
\providecommand{\MR}{\relax\ifhmode\unskip\space\fi MR }
\providecommand{\MRhref}[2]{%
  \href{http://www.ams.org/mathscinet-getitem?mr=#1}{#2}
}
\providecommand{\href}[2]{#2}

\end{document}